\documentclass[11pt]{amsart}

\usepackage{amssymb}
\usepackage{latexsym}
\usepackage{amsmath}

\begin{document}

\newtheorem{proposition}{Proposition}[section]
\newtheorem{definition}{Definition}[section]
\newtheorem{lemma}{Lemma}[section]
\newtheorem{theorem}{Theorem}[section]
\newtheorem{corollary}{Corollary}[section]
\newtheorem{problem}{Problem}[section]
\newtheorem{conjecture}{Conjecture}[section]

\title{Quantum groups and Hadamard matrices}
\author{Teodor Banica}
\address{T.B.: Department of Mathematics,
Paul Sabatier University, 118 route de Narbonne, 31062 Toulouse,
France} \email{banica@picard.ups-tlse.fr}
\author{Remus Nicoara}
\address{R.N.: Department of Mathematics, Vanderbilt University, 1326 Stevenson Center, Nashville, TN 37240, USA}
\email{remus.nicoara@vanderbilt.edu}
\thanks{R.N. was supported by NSF under Grant No. DMS-0500933}

\begin{abstract}
To any complex Hadamard matrix we associate a quantum permutation
group. The correspondence is not one-to-one, but the quantum group
encapsulates a number of subtle properties of the matrix. We
investigate various aspects of the construction: compatibility
to product operations, characterization of matrices which
give usual groups, explicit computations for small matrices.
\end{abstract}

\maketitle

\section*{Introduction}

A complex Hadamard matrix is a matrix $h\in M_n(\mathbb C)$,
having the following property: entries are on the unit circle, and
rows are mutually orthogonal. Equivalently, $n^{-1/2}h$ is a
unitary matrix with all entries having the same absolute value.

The basic example is the Fourier $n\times n$ matrix, given by
$F_{ij}=w^{ij}$, where $w=e^{2\pi i/n}$. There are many other
examples, but there is no other known family of complex Hadamard matrices that exist for every $n$.

These matrices are related to several questions:
\begin{enumerate}
\item  The story begins with work of Sylvester, who studied real
Hadamard matrices \cite{sy}. These have $\pm 1$ entries, and can
be represented as black and white pavements of a $n\times n$
square. Such a matrix must have order $n=2$ or $n=4k$, and the
main problem here is whether there are such matrices for any $k$.
At the time of writing, the verification goes up to $k=166$.

\item One can try to replace the $\pm 1$ entries by numbers of the
form $\pm 1,\pm i$, or by $n$-th roots of unity (the Fourier
matrix appears), or by roots of unity of arbitrary order. These
are natural generalizations of real Hadamard matrices, first
investigated by Butson in the sixties \cite{bu}.

\item The matrices with arbitrary complex entries appeared in
 the eighties. Popa discovered that such a matrix gives rise to the operator
algebra condition $\Delta\perp h\Delta h^*$, where $\Delta\subset
M_n(\mathbb C)$ is the algebra of diagonal matrices \cite{po1}.
Then work of Jones led to the conclusion that associated to $h$
are the following objects: a statistical mechanical model, a knot
invariant, a subfactor, and a planar algebra. The computation of
algebraic invariants for these objects appears to be remarkably
subtle and difficult. See \cite{j3}.

 \item In the
meantime, motivated by a problem of Enflo, Bj\"orck discovered
that circulant Hadamard matrices correspond to cyclic $n$-roots,
and started to construct many examples \cite{bj}. In a key paper,
Haagerup proved that for $n=5$ the Fourier matrix is the only
complex Hadamard matrix, up to permutations and multiplication by
diagonal unitaries \cite{h}. Self-adjoint complex Hadamard
matrices of order $6$ have been recently classified \cite{bn}.

\item Another important result is that of Petrescu, who discovered
a one-parameter family at $n=7$, providing a counterexample to a
conjecture of Popa regarding the finiteness of the number of
complex Hadamard matrices of prime dimension \cite{pe}. The notion
of deformation is further investigated in \cite{n}.

\item Several applications of complex Hadamard matrices were
discovered in the late nineties. These appear in a remarkable
number of areas of mathematics and physics, and the various
existence and classification problems are fast evolving. See
\cite{tz} for a catalogue of most known complex Hadamard matrices
and a summary of their applications.
\end{enumerate}

The purpose of this paper is to present a definition for what
might be called symmetry of a Hadamard matrix. This is done by
associating to $h$ a certain Hopf algebra $A$. This algebra is of
a very special type: it is a quotient of Wang's quantum
permutation algebra \cite{wa2}. In other words, we have the
heuristic formula $A=\mathbb C(G)$, where $G$ is a quantum group
which permutes the set $\{1,\ldots ,n\}$.

As an example, the Fourier matrix gives $G=\mathbb Z_n$. In the
general case $G$ doesn't really exist as a concrete object, but
some partial understanding of the quantum permutation phenomenon
is available, via methods from finite, compact and discrete
groups, subfactors, planar algebras, low-dimensional topology,
statistical mechanical models, classical and free probability,
geometry, random matrices etc. See \cite{bb}, \cite{bbc},
\cite{bc1}, \cite{bc2} and the references there. The hope is that
the series of papers on the subject will grow quickly, and evolve
towards the Hadamard matrix problematics. A precise comment in this sense is made in the last section.

Finally, let us mention that the fact that Hadamard matrices
produce quantum groups is known since \cite{ba1}, \cite{wa2}. The relation between coinvolutive compact quantum groups, abstract statistical mechanical models and associated commuting squares, subfactors and standard $\lambda$-lattices is worked out in \cite{ba2}. Another result is the one in \cite{ba3}, where a Tannaka-Galois type correspondence between spin planar algebras and quantum permutation groups is obtained. In both cases the situations discussed are more general than those involving Hadamard matrices, where some simplifications are expected to appear. However, the whole subject is quite technical, and it is beyond our purposes to discuss the global picture. This paper should be regarded as an introduction to the subject.

The paper is organized as follows. 1 is a quick introduction to
quantum permutation groups. 2 contains a few basic facts regarding
complex Hadamard matrices, the construction of the correspondence,
and the computation for the Fourier matrix. In 3-6 we investigate
several aspects of the correspondence: compatibility between
tensor products of quantum groups and of Hadamard matrices, the
relation with magic squares, characterization of matrices which
give usual groups, explicit computations for small matrices, and some
comments about subfactors and deformation.

\subsection*{Acknowledgements}
This work was started in March 2006 at Vanderbilt University, and
we would like to thank Dietmar Bisch for the kind hospitality and
help. The paper also benefited from several discussions with
Stefaan Vaes.

\section{Quantum permutation groups}

This section is an introduction to $A(S_n)$, the algebra of free
coordinates on the symmetric group $S_n$. This algebra was
discovered by Wang in \cite{wa2}.

The idea of noncommuting coordinates goes back to Heisenberg, the
specific idea of using algebras of free coordinates on algebraic
groups should be attributed to Brown, and a detailed study of
these algebras, from a K-theoretic point of view, is due to
McClanahan. Brown's algebras are in fact too big, for instance
they have no antipode, and the continuation of the story makes use
of Woronowicz's work on the axiomatization of compact
quantum groups. The first free quantum groups, corresponding to
$U_n$ and $O_n$, appeared in Wang's thesis. The specific question
about free analogues of $S_n$ was asked by Connes. See \cite{br},
\cite{m}, \cite{wa1}, \cite{wa2}.

We can see from this brief presentation that the algebra $A(S_n)$
we are interested in comes somehow straight from Heisenberg, after
a certain period of time. In fact, as we will see now, the
background required in order to define $A(S_n)$ basically reduces
to the definition of $C^*$-algebras, and to some early work on the
subject:

Let $A$ be a $C^*$-algebra. That is, we have a complex algebra
with a norm and an involution, such that Cauchy sequences
converge, and $||aa^*||=||a||^2$.

The basic example is $B(H)$, the algebra of bounded operators on a
Hilbert space $H$. In fact, any $C^*$-algebra appears as
subalgebra of some $B(H)$.

The key example is $C(X)$, the algebra of continuous functions on
a compact space $X$. This algebra is commutative, and any
commutative $C^*$-algebra is of this form.

The above two statements are results of Gelfand-Naimark-Segal and
Gelfand. The proofs make use of standard results in commutative
algebra, complex analysis, functional analysis, and measure
theory. The relation between the two statements is quite subtle,
and comes from the spectral theorem for self-adjoint operators.
The whole material can be found in any book on operator algebras.

Wang's idea makes a fundamental use of the notion of projection:

\begin{definition}
Let $A$ be a $C^*$-algebra.
\begin{enumerate}
\item A projection is an element $p\in A$ satisfying $p^2=p=p^*$.
\item Two projections $p,q\in A$ are called orthogonal when
$pq=0$. \item A partition of unity is a set of orthogonal
projections, which sum up to $1$.
\end{enumerate}
\end{definition}

A projection in $B(H)$ is an orthogonal projection $P(K)$, where
$K\subset H$ is a closed subspace. Orthogonality of projections
corresponds to orthogonality of subspaces, and partitions of unity
correspond to decompositions of $H$.

A projection in $C(X)$ is a characteristic function $\chi(Y)$,
where $Y\subset X$ is an open and closed subset. Orthogonality of
projections corresponds to disjointness of subsets, and partitions
of unity correspond to partitions of $X$.

\begin{definition}
A magic unitary is a square matrix $u\in M_n(A)$, all whose rows
and columns are partitions of unity in $A$.
\end{definition}

The terminology comes from a vague similarity with magic squares,
to be investigated later on. For the moment we are rather
interested in the continuing the above classical/quantum analogy, for projections and partitions of unity:

A magic unitary over $B(H)$ is of the form $P(K_{ij})$, with $K$
magic decomposition of $H$, in the sense that all rows and columns
of $K$ are decompositions of $H$. The basic examples here are of
the form $K_{ij}=\mathbb C\,\xi_{ij}$, where $\xi$ is a magic
basis of $H$, in the sense that all rows and columns of $\xi$ are
bases of $H$.

A magic unitary over $C(X)$ is of the form $\chi(Y_{ij})$, with
$Y$ magic partition of $X$, in the sense that all rows and columns
of $Y$ are partitions of $X$.  The key example here comes from a
finite group $G$ acting on a finite set $X$: the characteristic
functions $\chi_{ij}=\left\{\sigma\in G\mid \sigma(j)=i\right\}$
form a magic unitary over $C(G)$.

In the particular case of the symmetric group $S_n$ acting on
$\{1,\ldots ,n\}$, we have the following presentation result,
which follows from the Gelfand theorem:

\begin{theorem}
$C(S_n)$ is the universal commutative $C^*$-algebra generated by
$n^2$ elements $\chi_{ij}$, with relations making $(\chi_{ij})$ a
magic unitary matrix. Moreover, the maps
\begin{eqnarray*}
\Delta(\chi_{ij})&=&\sum \chi_{ik}\otimes \chi_{kj}\cr
\varepsilon(\chi_{ij})&=&\delta_{ij}\cr S(\chi_{ij})&=&\chi_{ji}
\end{eqnarray*}
are the comultiplication, counit and antipode of $C(S_n)$.
\end{theorem}

In other words, when regarding $S_n$ as an algebraic group, the
relations satisfied by the $n^2$ coordinates are those expressing
magic unitarity. Indeed, the characterstic functions $\chi_{ij}$ are nothing but the $n^2$ coordinates on the group $S_n\subset O_n$, where the embedding is the one given by permutation matrices.

See the preliminary section in \cite{bbc} for a
proof of the above result, and the introduction of \cite{bc2} for more comments on the subject.

We are interested in the algebra of free coordinates on $S_n$.
This is obtained by removing commutativity in the above
presentation result:

\begin{definition}
$A(S_n)$ is the universal $C^*$-algebra generated by $n^2$
elements $x_{ij}$, with relations making $(x_{ij})$ a magic
unitary matrix. The maps
\begin{eqnarray*}
\Delta(x_{ij})&=&\sum x_{ik}\otimes x_{kj}\cr
\varepsilon(x_{ij})&=&\delta_{ij}\cr S(x_{ij})&=&x_{ji}
\end{eqnarray*}
are called the comultiplication, counit and antipode of $A(S_n)$.
\end{definition}

This algebra, discovered by in \cite{wa2}, fits into the quantum
group formalism developed by Woronowicz in \cite{wo1}, \cite{wo2}. In fact, the quantum group $G_n$
defined by the formula $A(S_n)=C(G_n)$ is a free
analogue of the symmetric group $S_n$. This quantum group doesn't
exist of course: the idea is just that various properties of
$A(S_n)$ can be expressed in terms of it. As an example, the
canonical map $A(S_n)\to C(S_n)$ should be thought of as coming
from an embedding $S_n\subset G_n$.

\begin{proposition}
For $n=1,2,3$ we have $A(S_n)=C(S_n)$.
\end{proposition}

This follows from the fact that for $n\leq 3$, the entries of a
$n\times n$ magic unitary matrix have to commute. For instance at $n=2$ the matrix must be of the form
$$u=\begin{pmatrix}p&1-p\cr 1-p&p\end{pmatrix}$$
where $p$ is a projection, and entries of this matrix obviously commute.

The result is no longer true for $n=4$, where more complicated examples of magic unitary matrices are available, for instance via diagonal concatenation:
$$u=\begin{pmatrix}
p&1-p&0&0\cr
1-p&p&0&0\cr
0&0&q&1-q\cr
0&0&1-q&q
\end{pmatrix}$$

In fact, $A(S_n)$ with $n\geq 4$ is not commutative, and infinite
dimensional.

Consider now an arbitrary magic unitary matrix $u\in M_n(A)$. We say that $u$ is a
corepresentation of $A$ if its coefficients generate $A$, and if
the formulae
\begin{eqnarray*}
\Delta(u_{ij})&=&\sum u_{ik}\otimes u_{kj}\cr
\varepsilon(u_{ij})&=&\delta_{ij}\cr S(u_{ij})&=&u_{ji}
\end{eqnarray*}
define morphisms of $C^*$-algebras. These morphisms have by
definition values in the algebras $A\otimes A,\mathbb C$
 and $A^{op}$, and they are uniquely determined
 if they exist. In case they exist, these morphisms make $A$ into
 a Hopf algebra in the sense of \cite{wo1}.

This notion provides an axiomatization for quotients of $A(S_n)$:

\begin{definition}
A quantum permutation algebra is a $C^*$-algebra $A$, given with a
magic unitary corepresentation $u\in M_n(A)$.
\end{definition}

We have the following purely combinatorial approach to these algebras:

Recall first that the unitary representations of a discrete group
$L$ are in correspondence with representations of the group
algebra $C^*(L)$. The same is known to hold for discrete quantum
groups, so we have a correspondence between representations
\begin{eqnarray*}
L_n&\to& U(H)\cr A(S_n)&\to& B(H)\end{eqnarray*} where $L_n$ is
the discrete quantum group whose group algebra is $A(S_n)$. We can
therefore consider the quantum permutation algebra
$A=C^*\left(Im(L_n)\right)$.

All this is quite heuristic, but the construction of $A$ is
definitely possible:

\begin{definition}
Associated to a representation $\pi:A(S_n)\to B(H)$ is the
smallest quantum permutation algebra $A_\pi$ realizing a
factorization of $\pi$.
\end{definition}

In this definition we assume of course that the morphism making
$\pi$ factorize is the canonical one, given by $x_{ij}\to u_{ij}$.
The construction of $A_\pi$ is standard, for instance by taking
the quotient of $A(S_n)$ by an appropriate Hopf algebra ideal. The
uniqueness up to isomorphism is also clear. See \cite{ba1}.

We have the following examples:
\begin{enumerate}
\item For the counit representation $\pi:A(S_n)\to\mathbb C$ we
get $A_\pi=\mathbb C$. \item For a faithful representation
$\pi:A(S_n)\subset B(H)$ we get $A_\pi=A(S_n)$. \item If $A$ is a quantum permutation algebra, $A\subset B(H)$, and $\pi:A(S_n)\to A$ is a surjective representation of $A(S_n)$, then $A_\pi=A$.
\end{enumerate}

These examples are all trivial. Note however that the third one
has the following interesting consequence:

\begin{proposition}
Any quantum permutation algebra is of the form $A_\pi$.
\end{proposition}

\begin{proof}
Let $A$ be a quantum permutation algebra. We can compose the
canonical quotient map $A(S_n)\to A$ with an embedding $A\subset
B(H)$, say given by the GNS theorem, and we get a representation
$\pi$ as in the statement.
\end{proof}

\section{Hadamard matrices}

Let $h\in M_n(\mathbb C)$ be an Hadamard matrix. This means that all entries of $h$ have modulus $1$, and that rows of $h$ are mutually orthogonal. In other words, we have
$$h=\begin{pmatrix}h_1\cr h_2\cr \dots\cr h_n\end{pmatrix}$$
where the vectors $h_i$ are formed by complex numbers of modulus
$1$, and are orthogonal with respect to the usual scalar product
of $\mathbb C^n$, given by:
$$<x,y>=x_1\bar{y}_1+\ldots +x_n\bar{y}_n$$

It follows from definitions that the columns of $h$ are orthogonal
as well.

We have the following characterization of such matrices:

\begin{proposition}
For a matrix $h\in M_n(\mathbb C)$, the following are equivalent:
\begin{enumerate}
\item $h$ is an Hadamard matrix.
\item $n^{-1/2}h$ is a unitary matrix, all whose entries have same modulus.
\end{enumerate}
\end{proposition}

We should mention that in the operator algebra literature, it is rather $n^{-1/2}h$ that is called Hadamard. The other warning is that in the combinatorics literature, the Hadamard matrices are just the real ones, meaning those having $\pm 1$ entries. We hope that the hybrid terminology used in this paper won't cause any trouble.

We use capital letters to denote explicit Hadamard matrices. The
first example, which is in fact the only basic example, is the
Fourier matrix:

\begin{definition}
The Fourier matrix of order $n$ is given by $F_{ij}=w^{ij}$, where $w=e^{2\pi
i/n}$.
\end{definition}

We have a natural equivalence relation for Hadamard matrices,
given by $h\sim k$ if one can pass from $h$ to $k$ by permutations
of rows and columns, and by multiplications of rows and columns by
complex numbers of modulus $1$. See \cite{tz}.

For $n=1,2,3,5$ any Hadamard matrix is equivalent to the Fourier
one, see \cite{h}. At $n=4$ we have the following example, depending on $q$ on
the unit circle:
$$M_q=\begin{pmatrix}
1&1&1&1\cr 1&q&-1&-q\cr 1&-1&1&-1\cr 1&-q&-1&q
\end{pmatrix}$$

These are, up to equivalence, all $4\times 4$ Hadamard matrices.
As an example, the Fourier matrix corresponds to the value $q=\pm
i$. See Haagerup \cite{h}.

At $n=6$ we have several examples, as for instance:
$$H=\begin{pmatrix}
i&1&1&1&1&1\cr 1&i&1&-1&-1&1\cr 1&1&i&1&-1&-1\cr 1&-1&1&i&1&-1\cr
1&-1&-1&1&i&1\cr 1&1&-1&-1&1&i
\end{pmatrix}$$

This matrix appears in the paper of Haagerup \cite{h}. Another example appears in the paper of Tao \cite{t}. This is given by the following formula, with $w=e^{2\pi
i/3}$:
$$T=\begin{pmatrix}
1&1&1&1&1&1\cr 1&1&w&w&w^2&w^2\cr 1&w&1&w^2&w^2&w\cr
1&w&w^2&1&w&w^2\cr 1&w^2&w^2&w&1&w\cr 1&w^2&w&w^2&w&1
\end{pmatrix}$$

Observe that the Fourier, Haagerup and Tao matrices are based on
certain roots of unity. We have here the following notion:

\begin{definition}
The Butson class $H_l(n)$ consists of Hadamard matrices $h\in
M_n(\mathbb C)$ having the property $h_{ij}^l=1$ for any $i,j$.
\end{definition}

In this definition $l$ is a positive integer. It is convenient to
use as well the notation $H_\infty(n)$, for the class of all
$n\times n$ Hadamard matrices. Here are a few examples:
\begin{enumerate}
\item The Fourier matrix is in $H_n(n)$. \item The Haagerup matrix
is in $H_4(6)$. \item The Tao matrix is in $H_3(6)$. \item The
real Hadamard matrices are in $H_2(n)$. \item The matrix $M(q)$ is
in $H_l(4)$, where $l$ is the order of $q^2$.
\end{enumerate}

We call $n$ size, and $l$ level. In lack of some better idea, the
complexity of an Hadamard matrix will be measured by its size
and level, the size coming first.

At infinite level the self-adjoint matrices of order $6$ are
classified in \cite{bn}. The complete list of Hadamard matrices at
$n=6$ is not known.

We end this discussion with the following well-known fact:

\begin{proposition}
Each Hadamard matrix is equivalent to a Hadamard matrix which is dephased, meaning that the first row and column and the diagonal consist of $1$'s.
\end{proposition}

Consider now an arbitrary Hadamard matrix $h\in M_n(\mathbb C)$. The rows of $h$, denoted as usual $h_1,\ldots ,h_n$, can be regarded as elements of the algebra ${\mathbb C}^n$.

Since each $h_i$ is formed by complex numbers of modulus $1$, this element is invertible. We can therefore consider the following matrix of elements of $\mathbb C^n$:
$$\xi_{ij}=\frac{h_{j}}{h_i}$$

The scalar products on rows of $\xi$ are computed as follows:
\begin{eqnarray*}
<\xi_{ij},\xi_{ik}>&=&<h_j/h_i,h_k/h_i>\cr
&=&n<h_j,h_k>\cr &=&n\cdot\delta_{jk}
\end{eqnarray*}

In other words, each row of $\xi$ is an orthogonal basis of $\mathbb C^n$. A similar computation works for columns, so $\xi$ is a magic basis of $\mathbb C^n$. Thus we can apply the procedure from previous section, and we get a
magic unitary matrix, a representation, and a quantum permutation
algebra:

\begin{definition}
Let $h\in M_n(\mathbb C)$ be an Hadamard matrix. \begin{enumerate} \item
$\xi(h)$ is the magic basis given by
$\xi_{ij}=h_j/h_i$. \item $P(h)$ is the magic
unitary given by $P_{ij}=P(\xi_{ij})$. \item $\pi_h$ is the
representation given by $\pi_h(u_{ij})=P_{ij}$. \item $A(h)$ is the
quantum permutation algebra associated to $\pi$.
\end{enumerate}
\end{definition}

In other words, associated to $h$ are the rank one projections
$P(h_j/h_i)$, which form a magic unitary over $M_n(\mathbb C)$. We
consider the corresponding representation
$$\pi:A(S_n)\to M_n(\mathbb C)$$
we say that this comes from a representation $L_n\to U(n)$ of the
dual of the $n$-th quantum permutation group, and we consider the
algebra $A=C^*(Im(L_n))$.

It is routine to check that equivalent Hadamard matrices $h_1,
h_2$ give the same algebra, $A(h_1)=A(h_2)$. This is because the
standard equivalence operations, namely permutation of rows and
columns, and multiplication or rows and columns by scalars of
modulus $1$, give conjugate magic unitaries, hence conjugate
representations of $A(S_n)$. As for the converse, this is expected
not to hold, see the conclusion.

As a first example, consider the Fourier matrix:

\begin{theorem}
For the Fourier matrix $F$ of order $n$ we have $A(F)=C(\mathbb
Z_n)$.
\end{theorem}

\begin{proof}
We have the following computation, where $\rho=(w,w^2,\ldots
,w^n)$:
\begin{eqnarray*}
F_{ij}=w^{ij} &\Rightarrow& F_i=\rho^i\cr
&\Rightarrow&\xi_{ij}=\rho^{j-i}\cr&\Rightarrow&
P_{ij}=P(\rho^{j-i})
\end{eqnarray*}

Consider the cycle $c(i)=i-1$. We regard $\mathbb Z_n=\{0,1,\ldots
,n-1\}$ as a subgroup of $S_n$, via the embedding $k\to c^k$.
Consider the following diagram:
$$\begin{matrix}
A(S_n)&\ &\rightarrow&\ &M_n(\mathbb C)\cr \ \cr\downarrow&\ &\ &\
&\uparrow\cr \ \cr C(S_n)&\ & \longrightarrow&\ &C( \mathbb Z_n)
\end{matrix}$$

We define the connecting maps in the following way:
\begin{enumerate}
\item The map $A(S_n)\to M_n(\mathbb C)$ is the representation
$\pi$. \item The map $A(S_n)\to C(S_n)$ is the canonical one,
given by $x_{ij}\to\chi_{ij}$. \item The map $C(S_n)\to C(\mathbb
Z_n)$ is the transpose of $\mathbb Z_n\subset S_n$. \item The map
$C(\mathbb Z_n)\to M_n(\mathbb C)$ is given by $\delta_k\to
P(\rho^k)$.
\end{enumerate}

We compute the image of $\chi_{ij}$ by the third map:
\begin{eqnarray*}
{\rm Im}(\chi_{ij})&=&{\rm Im}\left(\chi\{\sigma\mid
\sigma(j)=i\}\right)\cr &=&\chi\{k\mid c^k(j)=i\}\cr
&=&\chi\{k\mid j-k=i\}\cr &=&\delta_{j-i}
\end{eqnarray*}

Thus at level of generators, we have the following diagram:
$$\begin{matrix}
x_{ij}&\ &\rightarrow&\ &P(\rho^{j-i})\cr \ \cr\downarrow&\ &\ &\
&\uparrow\cr \ \cr \chi_{ij}&\ & \longrightarrow&\ & \delta_{j-i}
\end{matrix}$$

This diagram commutes, so the above diagram of algebras commutes
as well. Thus we have a factorization of $\pi$ through the quantum
permutation algebra $C(\mathbb Z_n)$. Moreover, this algebra is
the minimal one making $\pi$ factorize, for instance because it is
isomorphic to the image of $\pi$. This gives the result.
\end{proof}

\section{Tensor products}

Given two Hadamard matrices $h\in M_n(\mathbb C)$ and $k\in
M_m(\mathbb C)$, we can consider their tensor product. This is an
element as follows:
$$h\otimes k\in M_n(\mathbb C)\otimes M_m(\mathbb C)$$

In order to view $h\otimes k$ as a usual matrix, we have to
identify the algebra on the right with the matrix algebra
$M_{nm}(\mathbb C)$. This isomorphism is not canonical, so we
proceed as follows: given $h,k$, we define $h\otimes k\in
M_{nm}(\mathbb C)$ by the formula
$$(h\otimes k)_{ia,jb}=h_{ij}k_{ab}$$
where double indices $ia,jb\in\{1,\ldots ,nm\}$ come from indices
$i,j\in \{1,\ldots ,n\}$ and $a,b\in\{1,\ldots m\}$, via some
fixed decomposition  of the index set. Most convenient here is to
use the lexicographic decomposition: $\{1,2,\ldots
,nm\}=\{11,12,\ldots ,nm\}$.

It is well-known that the above two definitions of $h\otimes k$
coincide.

\begin{proposition}
$h\otimes k$ is an Hadamard matrix.
\end{proposition}

\begin{proof}
We can use for instance Proposition 2.1: since $h,k$ are Hadamard
matrices, $n^{-1/2}h,m^{-1/2}k$ are unitaries. Now the tensor
product of two unitaries being a unitary, we get that
$(nm)^{-1/2}h\otimes k$ is unitary, so $h\otimes k$ is Hadamard.
\end{proof}

It is known from \cite{bb} that various operations at level of
combinatorial objects lead to similar operations at level of
quantum permutation algebras. In this section we prove such a
result for tensor products of Hadamard matrices:

\begin{theorem}
We have $A(h\otimes k)=A(h)\otimes A(k)$.
\end{theorem}

\begin{proof}
It is convenient here to use tensor product notations, at the same
time with bases and indices. We use the lexicographic
identification of Hilbert spaces
\begin{eqnarray*}\mathbb
C^n\otimes\mathbb C^m&=&\mathbb C^{nm}\cr e_i\otimes
e_a\,\,\,\,&=&e_{ia}
\end{eqnarray*}
along with the corresponding identification of operator algebras:
\begin{eqnarray*}
M_n(\mathbb C)\otimes M_m(\mathbb C)&=& M_{nm}(\mathbb C)\cr
e_{ij}\otimes e_{ab}\,\,\,\,\,\,\;\,\,\,\,&=&e_{ia,jb}
\end{eqnarray*}

We denote by $\xi,P,\pi$ the objects associated in Definition 2.3
to $h,k,h\otimes k$, with the choice of indices telling which is
which. We have the following computation:
\begin{enumerate}
\item  The $ia$-th row of $h\otimes k$ is given by $(h\otimes
k)_{ia}=h_i\otimes k_a$. \item The corresponding magic basis is
given by $\xi_{ia,jb}=\xi_{ij}\otimes\xi_{ab}$. \item The
corresponding magic unitary is given by $P_{ia,jb}=P_{ij}\otimes
P_{ab}$.
\end{enumerate}

Consider now the factorizations associated to $h,k$:
$$A(S_n)\to A(h)\to M_n(\mathbb C)$$
$$A(S_m)\to A(k)\to M_n(\mathbb C)$$

These are given by the following formulae:
$$x_{ij}\to u_{ij}\to P_{ij}$$
$$x_{ab}\to u_{ab}\to P_{ab}$$

We can form the tensor product of these maps, then we add to the
picture the factorization associated to $h\otimes k$, along with
two vertical arrows:
$$\begin{matrix}
A(S_n)\otimes A(S_m)&\to& A(h)\otimes A(k)&\to& M_n(\mathbb
C)\otimes M_m(\mathbb C)\cr \uparrow &&&&\downarrow\cr
A(S_{nm})&\to&A(h\otimes k)&\to&M_{nm}(\mathbb C)
\end{matrix}$$

Here the arrow on the right is the identification mentioned in the
beginning of the proof. As for the arrow on the left, this is
constructed as follows. Recall first that $A(S_{nm})$ is the
universal algebra generated by the entries of a $nm\times nm$
magic unitary matrix $x_{ia,jb}$. The matrix $x_{ij}\otimes
x_{ab}$ being a magic unitary as well, we get a morphism of
algebras mapping $x_{ia,jb}\to x_{ij}\otimes x_{ab}$, that we put
at left.

At level of generators, we have the following diagram:
$$\begin{matrix}
x_{ij}\otimes x_{ab}&\to&u_{ij}\otimes u_{ab}&\to& P_{ij}\otimes
P_{ab}\cr \uparrow &&&&\downarrow\cr
x_{ia,jb}&\to&u_{ia,jb}&\to&P_{ia,jb}
\end{matrix}$$

This diagram commutes, so the above diagram of algebras commutes
as well. Thus the factorization associated to $h\otimes k$
factorizes through the algebra $A(h)\otimes A(k)$. Now since
$A(h\otimes k)$ is the minimal algebra producing such a
factorization, this gives a morphism $A(h)\otimes A(k)\to
A(h\otimes k)$, mapping $u_{ij}\otimes u_{ab}\to u_{ia,jb}$.

This morphism is surjective, because the elements $u_{ia,jb}$
generate the algebra $A(h\otimes k)$. Now since $A(h\otimes k)$
makes a factorization of the representation associated to
$h\otimes k$, this algebra produces as well factorizations of the
representations associated to $h,k$. Thus our morphism is
injective, and this gives the result.
\end{proof}

As a first application, we can start classification work for
small Hadamard matrices. Recall from previous section that each
matrix has a size $n$ and a level $l$, according to the formula
$h\in H_l(n)$ with $l$ minimal, and that we decided to list
matrices in terms of $n,l$, with $n$ coming first. The situation
is as follows:
\begin{enumerate} \item For $n=1,2,3$ there is only one matrix,
namely the Fourier one. As shown by Theorem 2.1, this matrix
produces the algebra $C(\mathbb Z_n)$. \item For $n=4$ the
smallest possible level is $l=2$. We have here the matrices $M_q$,
with parameter $q=\pm 1$. \item The next possible level is $l=4$.
We have here the matrices $M_q$ with $q=\pm i$, both equivalent to
the Fourier matrix, which gives $C(\mathbb Z_4)$.
\end{enumerate}

Thus at $nl\leq 44$, we have just two matrices left. But these can
be investigated by using Theorem 3.1:

\begin{corollary}
For $q=\pm 1$, the Hadamard matrix
$$M_q=\begin{pmatrix} 1&1&1&1\cr 1&q&-1&-q\cr 1&-1&1&-1\cr
1&-q&-1&q
\end{pmatrix}$$
produces the algebra $C(\mathbb Z_2\times\mathbb Z_2)$.
\end{corollary}

\begin{proof}
The $2\times 2$ Fourier matrix is given by $F_{ij}=(-1)^{ij}$, so
we have:
$$F=\begin{pmatrix}-1&1\cr 1&1\end{pmatrix}$$

Consider the tensor square of $F$, given by $(F\otimes
F)_{ia,jb}=F_{ia}F_{jb}$. As for $F$, this is a $\pm 1$ matrix.
Now since the $-1$ entry of $F$ appears on the $11$ position, the
$-1$ entries of $F\otimes F$ appear on positions $(ia,jb)$
satisfying $ij=11$ or $ab=11$, with the position $(11,11)$
excluded. There are six such positions:
$$(11,12),(11,21),(12,11),(12,12),(21,11),(21,21)$$

Now by using the lexicographic identification
$\{1,2,3,4\}=\{11,12,21,22\}$, these six positions are
12,13,21,22,31,33. This gives the following formula:
$$F\otimes F=\begin{pmatrix} 1&-1&-1&1\cr -1&-1&1&1\cr -1&1&-1&1\cr
1&1&1&1
\end{pmatrix}$$

By permuting rows and columns of this matrix we can get both
matrices in the statement. The result follows from Theorem 2.1 and
Theorem 3.1, by using the canonical identification $C(\mathbb
Z_2)\otimes C(\mathbb Z_2)=C(\mathbb Z_2\times\mathbb Z_2)$.
\end{proof}

\section{Magic squares}

In this section and in next one we investigate the case when $A$ is commutative.
In this situation, we must have $A=C(G)$, for a certain subgroup
$G\subset S_n$.

We already know that $A$ is commutative for the Fourier matrix,
where we have $G=\mathbb Z_n$, and for the matrix $M_q$ with
$q=\pm 1$, where we have $G=\mathbb Z_2\times\mathbb Z_2$. In
fact, these two examples are the only ones that have been
investigated so far. This makes our general problem quite unclear,
because we have no counterexample.

Let us go back however to proof of Theorem 2.1. A careful
examination shows that the proof relies on the fact that $P$ is a
circulant matrix. The idea will be to generalize this proof, with
``circulant'' replaced by ``magic''.

Let us first recall the following well-known definition:

\begin{definition}
A magic square is a square matrix $\sigma$, having as entries the
numbers $1,\ldots ,n$, such that each row and each column of
$\sigma$ is a permutation of $1,\ldots ,n$.

We assume that all magic squares are normalized, in the sense that
the first row is $(1,\ldots ,n)$, and the diagonal is $(1,\ldots
,1)$.
\end{definition}

In this definition, we use the fact that it is always possible to
normalize a magic square: a permutation of the columns makes the
first row $(1,\ldots ,n)$, then a permutation of the rows makes
the diagonal $(1,\ldots ,1)$.

As a first example, we have the $n\times n$ circulant matrix
$\sigma(i,j)=j-i$, with $j-i\in\{1,\ldots,n\}$ taken modulo $n$.
Here is another example:
$$\sigma=\begin{pmatrix}1&2&3&4\cr 2&1&4&3\cr 3&4&1&2\cr
4&3&2&1\end{pmatrix}$$

The magic squares produce magic unitaries, in the following way:

\begin{proposition}
If $E=(E_1,\ldots ,E_n)$ is a partition of unity with rank one
projections and $\sigma=\sigma(i,j)$ is a magic square, then
$(E_\sigma)_{ij}=E_{\sigma(i,j)}$ is a magic unitary.
\end{proposition}

In other words, $E_\sigma$ is obtained by putting $E$ superscripts
to elements of $\sigma$. For instance the above $4\times 4$ matrix
gives:
$$E_\sigma=\begin{pmatrix}E_1&E_2&E_3&E_4\cr E_2&E_1&E_4&E_3\cr
 E_3&E_4&E_1&E_2\cr
E_4&E_3&E_2&E_1\end{pmatrix}$$

In the above statement, the condition that $\sigma$ is normalized
is not necessary; nor the fact that the projections $E_i$ are of
rank one. We make these assumptions for some technical reasons, to
become clear later on.

We denote by $\sigma_1,\ldots ,\sigma_n$ the rows of a magic
square $\sigma$, and we regard them as permutations of $\{1,\ldots
,n\}$. For instance, for the above $4\times 4$ matrix we get:
\begin{eqnarray*}
\sigma_1&=&\begin{pmatrix}1&2&3&4\end{pmatrix}\cr
\sigma_2&=&\begin{pmatrix}2&1&4&3\end{pmatrix}\cr
\sigma_3&=&\begin{pmatrix}3&4&1&2\end{pmatrix}\cr
\sigma_4&=&\begin{pmatrix}4&3&2&1\end{pmatrix}
\end{eqnarray*}

With these notations, we have the following result about
commutativity, which deals with a slightly more general situation
than that of Hadamard matrices:

\begin{theorem}
For a representation $\pi:A(S_n)\to M_n(\mathbb C)$ having the
property that $P_{ij}=\pi(x_{ij})$ are rank one projections, the
following are equivalent:
\begin{enumerate}
\item $A_\pi$ is commutative. \item $A_\pi=C(G)$, for a certain
subgroup $G\subset S_n$. \item $P=E_\sigma$, for a certain
partition of unity $E$ and magic matrix $\sigma$.
\end{enumerate}

Moreover, in this situation $G$ is the group generated by the rows
of $\sigma$.
\end{theorem}

\begin{proof}
Assume that $A=A_\pi$ is commutative. The algebra ${\rm Im}(\pi)$
being a quotient of $A$, it is commutative as well. Thus the
projections $P_{ij}$ mutually commute.

On the other hand, two rank one projections commute if and only if
their images are equal, or orthogonal. Together with the magic
unitarity of $P$, this shows that each row of $P$ is a permutation
of the first row.

So, consider the first row of $P$, regarded as a partition of
unity:
$$E=(E_1,\ldots ,E_n)$$

By the above remark, the $i$-th row of $P$ must be of the
following form:
$$E_{\sigma_i}=(E_{\sigma_i(1)},\ldots ,E_{\sigma_i(n)})$$

Here $\sigma_i\in S_n$ is a certain permutation. Now the formula
$\sigma(i,j)=\sigma_i(j)$ defines a matrix $\sigma$, which is
magic. With the above $E$ and this matrix $\sigma$, we have
$P=E_\sigma$.

In other words, we have (1) $\Rightarrow$ (3). Also, the last
assertion implies (2), which in turn implies (1). So, it remains
to prove that (3) implies the last assertion.

Consider the following diagram, where
$X=\{\sigma_1,\ldots,\sigma_n\}$ is the subset of $S_n$ formed by
rows of $\sigma$:
$$\begin{matrix}
A(S_n)&\ &\rightarrow&\ &M_n(\mathbb C)\cr \ \cr\downarrow&\ &\ &\
&\uparrow\cr \ \cr C(S_n)&\ & \longrightarrow&\ &C(X)
\end{matrix}$$

We define the connecting maps in the following way:
\begin{enumerate}
\item The map $A(S_n)\to M_n(\mathbb C)$ is the representation
$\pi$. \item The map $A(S_n)\to C(S_n)$ is the canonical one,
given by $x_{ij}\to\chi_{ij}$. \item The map $C(S_n)\to C(X)$ is
the transpose of $X\subset S_n$. \item The map $C(X)\to
M_n(\mathbb C)$ is given by $\delta_{\sigma_k}\to E_k$.
\end{enumerate}

We compute the image of $\chi_{ij}$ by the third map:
\begin{eqnarray*}
{\rm Im}(\chi_{ij})&=&{\rm Im}\left(\chi\{\sigma\mid
\sigma(j)=i\}\right)\cr &=&\chi\{\sigma_k\mid \sigma_k(j)=i\}\cr
&=&\delta_{\sigma_{\sigma(i,j)}}
\end{eqnarray*}

Thus at level of generators, we have the following diagram:
$$\begin{matrix}
x_{ij}&\ &\rightarrow&\ &E_{\sigma(i,j)}\cr \ \cr\downarrow&\ &\
&\ &\uparrow\cr \ \cr \chi_{ij}&\ & \longrightarrow&\ &
\delta_{\sigma_{\sigma(i,j)}}
\end{matrix}$$

This diagram commutes, so the above diagram of algebras commutes
as well. Now since $C(S_n)$ is a Hopf algebra, the Hopf algebra
$A_\pi$ we are looking for must be a quotient of it. In other
words, we must have $A=C(G)$, where $G$ is a subgroup of $S_n$. On
the other hand, from minimality of $A_\pi$ we get that this
algebra must be the minimal one containing $C(X)$. Thus $G$ is the
group generated by $X$.
\end{proof}

As an illustration, we get new proofs for Theorem 2.1 and
Corollary 3.1:

\begin{enumerate}
\item For the Fourier matrix we have $P=E_\sigma$, where $\sigma$
is the magic square given by $\sigma(i,j)=j-i$, and $E$ is the
first row of $P$. We have $\sigma_i=c^i$, where $c$ is the cycle
$c(i)=i-1$, so we get the group $\{1,c,\ldots ,c^{n-1}\}=\mathbb
Z_n$. \item For the matrix $M_q$ with $q=\pm 1$ we have
$P=E_\sigma$, where $\sigma$ is the $4\times 4$ matrix in the
beginning of this section, and $E$ is the first row of $P$. We get
the group formed by rows of $\sigma$, namely
$\{\sigma_1,\sigma_2,\sigma_3,\sigma_4\}\simeq \mathbb
Z_2\times\mathbb Z_2$.
\end{enumerate}

In the general case, the construction $\sigma\to\pi\to A_\pi$ is
not understood. Some comments in this sense are presented in the
end of next section.

\section{The commutative case}

We have now all ingredients for answering the question raised in
the beginning of previous section. The result is as follows:

\begin{theorem}
For an Hadamard matrix $h$, the following are equivalent:
\begin{enumerate}\item $A$ is commutative. \item $h$ is a tensor
product of Fourier matrices.\end{enumerate}
\end{theorem}

\begin{proof}
The statement is of course up to equivalence of Hadamard matrices.
We use this fact several times in the proof, without special
mention: that is, we agree to allow suitable permutations of rows
and columns of $h$, as well as multiplication of rows and columns
by scalars of modulus one.

(2) $\Rightarrow$ (1): this follows from Theorem 2.1 and Theorem
3.1.

(1) $\Rightarrow$ (2): assume that $A$ is commutative. Theorem 4.1
tells us that we have $P=E_\sigma$, for a certain magic square
$\sigma$. Here $E$ is a certain partition of unity, but since
$\sigma$ is normalized, $E$ must be the first row of $P$.

We have $P_{ij}=P(h_j/h_i)$, so the condition $P=E_\sigma$
becomes:
$$P\left(\frac{h_j}{h_i}\right)=P\left(\frac{h_{\sigma(i,j)}}{h_1}\right)$$

By using Proposition 2.2. we can assume $h_1=1$. This gives:
$$P\left(\frac{h_j}{h_i}\right)=P\left(h_{\sigma(i,j)}\right)$$

Now two rank one projections being equal if and only if the
corresponding vectors are proportional, we must have complex
scalars $\lambda_{ij}$ such that:
$$\frac{h_j}{h_i}=\lambda_{ij}\,h_{\sigma(i,j)}$$

Now remember that the row vectors $h_i$ have entries of modulus
one, so they are elements of the group $\mathbb T^n$, where
$\mathbb T$ is the unit circle. We denote by $H_i\in \mathbb
T^n/\mathbb T$ their images modulo constant vectors in $\mathbb
T^n$. The above condition becomes:
$$\frac{H_j}{H_i}=H_{\sigma(i,j)}$$

This relation shows that $X=\{H_1,\ldots ,H_n\}$ is stable by
quotients. Since $X$ contains as well the neutral element $H_1$,
it must be a subgroup of $\mathbb T^n/\mathbb T$.

Now since $X$ is an abelian group, it must be a product of cyclic
groups, and we can proceed as follows:

(1) Assume first that we have $X\simeq\mathbb Z_n$. We can permute
rows of $h$, as to have $H_i=H^{i-1}$ for any $i$, for a certain
element $H\in \mathbb T^n/\mathbb T$. This gives
$h_i=\lambda_i\rho^{i-1}$ for any $i$, for a certain element
$\rho\in\mathbb T^n$, and certain scalars $\lambda_i$. Now since
the vectors $h_i$ are orthogonal, their multiples $\rho^{i-1}$
must be orthogonal as well.

Consider now the vector $\rho$. Its image in $\mathbb T^n/\mathbb
T$ has order $n$, so by multiplying $h$ by a suitable scalar, we
may assume that we have $\rho^n=1$. In other words, we have
$\rho=(w_i)$, where each $w_i$ is a $n$-th root of unity.

Now from the relation $\rho^{i-1}\perp \rho^{j-1}$ for $1\leq
i<j\leq n$ we get that the sum of $k$-th powers of the numbers
$w_i$ vanishes, for $1\leq k<n$. This shows that the product of
degree one polynomials $X-w_i$ is the polynomial $X^n-1$, so the
set of numbers $w_i$ is the set of $n$-th roots of unity. Thus by
permuting rows of $h$ we can assume that we have $w_i=w^i$, with
$w=e^{2\pi i/n}$.

Now the relation $h_i=\lambda_i\rho^{i-1}$ tells us that $h$ is
obtained from the Fourier matrix by multiplying rows by scalars.
Thus $h$ is equivalent to the Fourier matrix.

(2) Assume now that we have $X\simeq Y\times Z$, for certain
groups $Y,Z$. By replacing $n$ by $nm$ and single indices by
double indices, we can assume that we have $H_{ia}=K_iL_a$, where
$Y=\{K_1,\ldots ,K_n\}$ and $Z=\{L_1,\ldots ,L_m\}$ are subgroups
of $\mathbb T^{nm}/\mathbb T$.

We take now arbitrary lifts $k_1,\ldots ,k_n$ and $l_1,\ldots
,l_m$ for the elements of $Y,Z$. The formula $H_{ia}=K_iL_a$
becomes $h_{ia}=\lambda_{ia}k_il_a$, for certain scalars
$\lambda_{ia}$. Now these vectors being orthogonal, by keeping $a$
fixed we get that the vectors $k_i$ are orthogonal, so $k$ is an
Hadamard matrix. Also, with $i$ fixed we get that $l$ is an
Hadamard matrix.

On the other hand, we can get rid of the scalars $\lambda_{ij}$ by
multiplying rows of $h$ by their inverses. Thus we have
$h_{ia}=k_il_a$, where $k,l$ are Hadamard matrices.

With notations from previous section, this gives $h=k\otimes l$.

(3) We can conclude now by induction. Assume that the statement is
proved for $n<N$, and let $h$ be a Hadamard matrix of order $N$.
In the case $X\simeq\mathbb Z_N$ we can apply (1) and we are done.
If not, we have $X\simeq Y\times Z$ as in (2), so we get
$h=k\otimes l$. Now Theorem 3.1 gives $A(h)=A(k)\otimes A(l)$, so
both $A(k),A(l)$ are commutative. Thus we can apply the induction
assumption, and we are done.
\end{proof}

The permutation groups coming from arbitrary magic squares can be
non-abelian. Consider for instance the following matrix:
$$\sigma=\begin{pmatrix}1&2&3&4&5\cr 3&1&2&5&4\cr 4&5&1&3&2\cr
2&4&5&1&3\cr 5&3&4&2&1\end{pmatrix}$$

The rows of this matrix generate the group $S_5$, which cannot be
obtained by using $5\times 5$ Hadamard matrices, because of
\cite{h}. In general, we don't know what are the groups which can
appear from the magic square construction.

The other question is about what happens to Theorem 4.1 when
the rank one assumption is removed. Once again, we don't know the
answer.

\section{Small matrices}

We recall from \cite{h} that for $n=1,2,3,5$ the only Hadamard
matrix is the Fourier one, and that at $n=4$ we have only the
matrices $M_q$, with $|q|=1$:
$$M_q=\begin{pmatrix}
1&1&1&1\cr 1&q&-1&-q\cr 1&-1&1&-1\cr 1&-q&-1&q
\end{pmatrix}$$

By results in previous sections, for parameters $q$ satisfying
$q^4=1$, this matrix $M_q$ produces a commutative algebra:
\begin{enumerate} \item  For $q=\pm 1$ we get $C(G)$, with $G=\mathbb
Z_2\times\mathbb Z_2$. \item For $q=\pm i$ we get $C(G)$, with
$G=\mathbb Z_4$.
\end{enumerate}

Now remember that for a finite abelian group $G$, the algebra of
complex functions $C(G)$ is canonically isomorphic to the group
algebra $C^*(G)$. The isomorphism is given by the discrete Fourier
transform in the case $G=\mathbb Z_n$, and by a product of such
transforms, or just by Pontrjagin duality, in the general case.

We can therefore reformulate the above result in the following
way:
\begin{enumerate} \item  For $q=\pm 1$ we get $C^*(G)$, with $G=\mathbb
Z_2\times\mathbb Z_2$. \item For $q=\pm i$ we get $C^*(G)$, with
$G=\mathbb Z_4$.
\end{enumerate}

In order to be generalized, this result has to reformulated one
more time. We have the following equivalent statement, in terms of
crossed products:
\begin{enumerate} \item  For $q=\pm 1$ we get $C^*(G)$, with $G=\mathbb
Z_2\rtimes\mathbb Z_2$. \item For $q=\pm i$ we get $C^*(G)$, with
$G=\mathbb Z_1\rtimes\mathbb Z_4$.
\end{enumerate}

Recall now that for a discrete group $G$, meaning a possibly
infinite group, without topology on it, the group algebra $C^*(G)$
is a $C^*$-algebra, obtained from the usual group algebra $\mathbb
C[G]$ by using a standard completion procedure.

With these preliminaries and notations, we have the following
result:

\begin{theorem}
Let $n$ be the order of $q^2$, and for $n<\infty$ write $n=2^sm$,
with $m$ odd. The matrix $M_q$ produces the algebra $C^*(G)$,
where $G$ is as follows:
\begin{enumerate}
\item For $s=0$ we have $G=\mathbb Z_{2n}\rtimes\mathbb Z_2$.
\item For $s=1$ we have $G=\mathbb Z_{n/2}\rtimes\mathbb
Z_4$.\item For $s\geq 2$ we have $G=\mathbb Z_{n}\rtimes\mathbb
Z_4$.\item For $n=\infty$ we have $G=\mathbb Z\rtimes\mathbb Z_2$.
\end{enumerate}
\end{theorem}

\begin{proof}
Consider the second and the third row of $M_q$:
\begin{eqnarray*}
\rho&=&\begin{pmatrix}1&q&-1&-q\end{pmatrix}\cr
\varepsilon&=&\begin{pmatrix}1&-1&1&-1\end{pmatrix}
\end{eqnarray*}

We regard these vectors as elements of the algebra $\mathbb C^4$.
The product of these vectors being the last row of $M_q$, we have
the following formula:
$$M_q=\begin{pmatrix}1\cr\rho\cr\varepsilon\cr
\varepsilon\rho\end{pmatrix}$$

By taking quotients of these vectors we get the corresponding
magic basis:
$$\xi=\begin{pmatrix}
1&\rho&\varepsilon&\varepsilon\rho\cr
\rho^{-1}&1&\varepsilon\rho^{-1}&\varepsilon\cr
\varepsilon&\varepsilon\rho&1&\rho\cr\varepsilon\rho^{-1}&
\varepsilon&\rho^{-1}&1
\end{pmatrix}$$

Each entry of this matrix is of the form $\rho^k$ or
$\varepsilon\rho^k$, with $k\in\mathbb Z$. So, consider the
vectors $\rho^k_+=\rho^k$ and $\rho^k_-=\varepsilon\rho^k$. These
are given by the following global formula:
$$\rho^k_\pm=\begin{pmatrix}1& \pm q^k&(-1)^k&\pm
(-q)^k\end{pmatrix}$$

Now the orthogonal projection onto a vector $(x_i)$ being the
matrix $(\bar{x}_ix_j)$, we have the following formula for the
corresponding projection:
$$P(\rho^k_\pm)=\begin{pmatrix}1&\pm q^{-k}&(-1)^k
&\pm (-q)^{-k}\cr \pm q^k&1&\pm (-q)^k&(-1)^k\cr (-1)^k&\pm
(-q)^{-k}&1&\pm q^{-k}\cr \pm (-q)^k&(-1)^k&\pm
q^k&1\end{pmatrix}$$

The idea is to express this matrix as a suitable linear
combination of elements of $U(4)$. We write $n=2^sm$ as in the
statement, and we proceed as follows:

\smallskip

\underline{Case $s=0$}. Here $n$ is odd, and we have $q^{2n}=1$,
so $q^n=\pm 1$. We can replace if necessary $q$ by $-q$, as to
have $q^n=1$.

Consider the canonical subgroup $\mathbb Z_2\times\mathbb
Z_2\subset U(4)$, consisting of the identity matrix $1$, and of
the following three matrices:
$$\alpha=\begin{pmatrix}0&1&0&0\cr 1&0&0&0\cr 0&0&0&1\cr
0&0&1&0\end{pmatrix}\hskip 5mm\beta=\begin{pmatrix}0&0&1&0\cr
0&0&0&1\cr 1&0&0&0\cr 0&1&0&0\end{pmatrix}\hskip
5mm\gamma=\begin{pmatrix}0&0&0&1\cr 0&0&1&0\cr 0&1&0&0\cr
1&0&0&0\end{pmatrix}$$

Consider also the following matrix, which is in $U(4)$ as well:
$$\sigma=\begin{pmatrix}0&0&q&0\cr
0&0&0&q^{-1}\cr q&0&0&0\cr 0&q^{-1}&0&0\end{pmatrix}$$

We have the following formulae, valid for $k$ odd:
$$\sigma^k=\begin{pmatrix}0&0&q^k&0\cr 0&0&0&q^{-k}\cr
q^k&0&0&0\cr 0&q^{-k}&0&0\end{pmatrix}\hskip
5mm\alpha\sigma^k=\begin{pmatrix}0&0&0&q^{-k}\cr 0&0&q^k&0\cr
0&q^{-k}&0&0\cr q^k&0&0&0\end{pmatrix}$$

We also have the following formulae, valid for $k$ even:
$$\sigma^k=\begin{pmatrix}q^k&0&0&0\cr 0&q^{-k}&0&0\cr
0&0&q^k&0\cr 0&0&0&q^{-k}\end{pmatrix}\hskip
5mm\alpha\sigma^k=\begin{pmatrix}0&q^{-k}&0&0\cr q^k&0&0&0\cr
0&0&0&q^{-k}\cr 0&0&q^k&0\end{pmatrix}$$

In particular we have $\sigma^n=\beta$ and $\sigma^{2n}=1$.

For $k$ odd the number $n+k$ is even, and we have:
\begin{eqnarray*}
P(\rho_\pm^k)&=&1\pm\alpha\sigma^{n+k}-\sigma^n\mp
\alpha\sigma^k\cr &=& (1-\sigma^n)(1\mp
\alpha\sigma^k)\end{eqnarray*}

For $k$ even the number $n+k$ is odd, and we have:
\begin{eqnarray*}
P(\rho_\pm^k)&=&1\pm\alpha\sigma^{k}+\sigma^n\mp
\alpha\sigma^{n+k}\cr &=& (1+\sigma^n)(1\pm
\alpha\sigma^k)\end{eqnarray*}

This gives the following global formula:
$$P(\rho_\pm^k)=(1+(-1)^k\sigma^n)(1\pm
(-1)^k\alpha\sigma^k)$$

It is routine to check that $\alpha,\sigma\in U(4)$ generate the
dihedral group $G=\mathbb Z_{2n}\rtimes\mathbb Z_2$. Consider now
the following element of the abstract group algebra $C^*(G)$:
$$P^k_\pm=
(1+(-1)^k\sigma^n)(1\pm(-1)^k\alpha\sigma^k)$$

From $\sigma^{2n}=1$ we get that $\sigma^n$ is self-adjoint, and
of square $1$. From $\alpha^2=1$ and $(\alpha\sigma)^2=1$ we get
$\sigma\alpha\sigma=\alpha$, hence $\sigma^k\alpha\sigma^k=\alpha$
for any $k$. It follows that $\alpha\sigma^k$ is also
self-adjoint, and of square $1$. This shows that in the above
formula, both elements in brackets are projections. Moreover,
these two projections commute. It follows that each element
$P_\pm^k$ is a projection.

Now remember that the magic unitary associated to $\pi$ has as
entries elements of the form $P(\rho^k_\pm)$. By making the
replacement $P(\rho^k_\pm)\to P^k_\pm$ we get a square matrix over
$C^*(G)$, all whose entries are projections. The sums on rows and
columns being $1$, this is a magic unitary, and we get a
factorization of $\pi$ through $C^*(G)$.

\smallskip

\underline{Case $n=\infty$}. What changes in the above proof is
that $\alpha,\sigma\in U(4)$ generate now the infinite dihedral
group $G=\mathbb Z\rtimes\mathbb Z_2$. The second part of the
proof applies as well, with this modification, and gives the
result.

\smallskip

\underline{Case $s=1$}. Here $n/2$ is odd. We have $q^{n}=-1$, so
$q^{n/2}=\pm i$, which gives $(-iq)^{n/2}=\pm 1$. We replace if
necessary $q$ by $-q$, as to have $(-iq)^{n/2}=1$.

Consider the canonical subgroup $\mathbb Z_4\subset U(4)$,
consisting of the identity matrix $1$, and of the following three
matrices:
$$\delta=\begin{pmatrix}0&1&0&0\cr 0&0&1&0\cr 0&0&0&1\cr
1&0&0&0\end{pmatrix}\hskip 5mm\delta^2=\begin{pmatrix}0&0&1&0\cr
0&0&0&1\cr 1&0&0&0\cr 0&1&0&0\end{pmatrix}\hskip
5mm\delta^3=\begin{pmatrix}0&0&0&1\cr 1&0&0&0\cr 0&1&0&0\cr
0&0&1&0\end{pmatrix}$$

Consider also the following matrix, which is in $U(4)$ as well:
$$\tau=\begin{pmatrix}-q&0&0&0\cr 0&q^{-1}&0&0\cr 0&0&-q&0\cr
0&0&0&q^{-1}\end{pmatrix}$$

Since $n/2$ is odd and $(-iq)^{n/2}=1$, we have $(i\tau)^{n/2}=1$.

We have the following formulae:
$$\delta\tau^k=\begin{pmatrix}0&q^{-k}&0&0\cr 0&0&(-q)^{k}&0\cr
0&0&0&q^{-k}\cr (-q)^k&0&0&0\end{pmatrix}$$
$$\delta^3\tau^k=\begin{pmatrix}0&0&0&q^{-k}\cr (-q)^{k}&0&0&0\cr
0&q^{-k}&0&0\cr 0&0&(-q)^k&0\end{pmatrix}$$

This gives the following formula for the above projection:
\begin{eqnarray*}
P(\rho_\pm^k)&=&1\pm\delta\tau^k+(-1)^k\delta^2\pm
(-1)^k\delta^3\tau^k\cr &=& (1+(-1)^k\delta^2)(1\pm
\delta\tau^k)\end{eqnarray*}

Consider now the matrix $\nu=i\tau$. We have $\nu^{n/2}=1$, and it
is routine to check that $\delta,\nu\in U(4)$ generate the group
$G=\mathbb Z_{n/2}\rtimes\mathbb Z_4$. Consider now the following
element of the abstract group algebra $C^*(G)$, where
$\tau=-i\nu$:
\begin{eqnarray*}
P^k_\pm&=&(1+(-1)^k\delta^2)(1\pm\delta\tau^k)
\end{eqnarray*}

From $\delta^4=1$ and from $\tau\delta\tau=-\delta$ we get that
each $P_\pm^k$ is a projection.

Now remember that the magic unitary associated to $\pi$ has as
entries elements of the form $P(\rho^k_\pm)$. By making the
replacement $P(\rho^k_\pm)\to P^k_\pm$ we get a square matrix over
$C^*(G)$, all whose entries are projections. The sums on rows and
columns being $1$, this is a magic unitary, and we get a
factorization of $\pi$ through $C^*(G)$.

\smallskip

\underline{Case $s\geq 2$}. Here $n/2$ is even, and we have
$q^{n}=-1$, so $q^{n/2}=\pm i$.

We use the matrices $\delta,\tau$. The above formulae for
$\delta\tau^k$ and $\delta^3\tau^k$ hold again, and lead to the
above formula for $P(\rho_\pm^k)$. What changes is the $n/2$-th
power of $\tau$:
$$\tau^{n/2}=\begin{pmatrix}\pm i&0&0&0\cr 0&\mp i&0&0\cr 0&0&\pm i&0\cr
0&0&0&\mp i\end{pmatrix}$$

Observe also that we have $\tau^n=-1$, so $(w\tau)^n=1$, where
$w=e^{\pi i/n}$.

This shows that we don't have a factorization like in the case
where $n/2$ is odd, so we must use the group generated by
$\delta,w\tau\in U(4)$, which is $G=\mathbb Z_n\rtimes\mathbb
Z_4$.
\end{proof}

The above result has several interpretations. First, it provides examples of subgroups of the $4$-th quantum permutation group, also called Pauli quantum group \cite{bc2}. In other words, the whole thing should be regarded as being part of a general McKay correspondence for the Pauli quantum group, not available yet.

Another interpretation is in terms of subfactors: the above result can be regarded as being a particular case of the $ADE$ classification of index $4$ subfactors. Unlike in the quantum group case, the classification is available here in full generality. See \cite{ek}. The principal graphs corresponding to the matrices $M_q$ are those of type $A_4^{(1)}$, $D_n^{(1)}$, $D_\infty$. This follows for instance by combining Theorem 6.1 with general results in \cite{ba2}, and with the well-known fact that cyclic groups correspond to type $A$ graphs, and dihedral groups correspond to type $D$ graphs.

As a last remark, the above result provides the first example of a deformation
situation for quantum permutation groups. This is of course just an example, but the general fact that it suggests would be use of the unit circle as parameter space, for more general deformation situations. We should mention that this idea,
while being fundamental in most theories emerging from Drinfeld's
original work \cite{d}, is quite new in the compact quantum group
area, where the deformation parameter traditionally belongs to the
real line. It is our hope that further developments of the subject
will be of use in connection with several problems, regarding both
Hadamard matrices and compact quantum groups.

\end{document}